\documentclass[final,12pt,3p]{elsarticle}

\usepackage[T1]{fontenc}
\usepackage[utf8]{inputenc}
\usepackage{amsmath,amssymb,amsthm,mathtools}
\usepackage{booktabs}
\usepackage{graphicx}
\usepackage{caption}

\usepackage[shortlabels]{enumitem}
\usepackage{geometry}
\usepackage{lineno}

\usepackage{xcolor}
\usepackage[expansion=false]{microtype}
\usepackage{hyperref} 

\biboptions{sort&compress}
\hypersetup{
  colorlinks=true,
  linkcolor=blue,
  citecolor=blue,
  urlcolor=blue,
  pdftitle={Unconditionally Positive Nonstandard Finite Difference Schemes for p-Laplacian Reaction--Diffusion Equations},
  pdfauthor={Achraf Zinihi, Matthias Ehrhardt, Moulay Rchid Sidi Ammi},
}
\allowdisplaybreaks 

\theoremstyle{plain}
\newtheorem{theorem}{Theorem}[section]
\newtheorem{proposition}[theorem]{Proposition}

\theoremstyle{remark}
\newtheorem{remark}[theorem]{Remark}

\newcommand{\eps}{\varepsilon}

\journal{}
\begin{document}

\begin{flushright}
		{\includegraphics[scale=0.2]{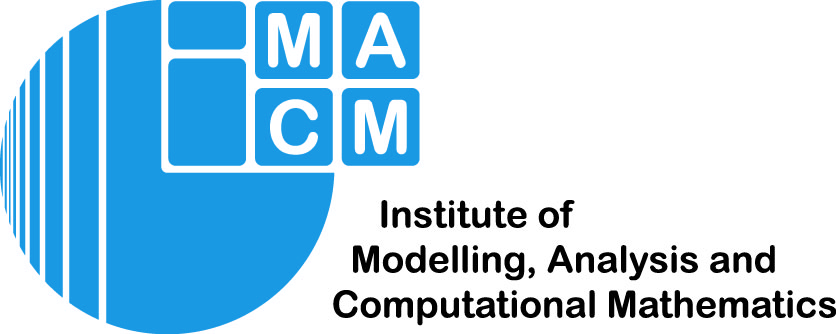}}
\end{flushright}

\begin{frontmatter}

\title{Unconditionally Positive Nonstandard Finite Difference Schemes\\  for $p$-Laplacian Reaction--Diffusion Equations}

\author[BUW,UMI]{Achraf Zinihi}
\ead{achrafzinihi99@gmail.com}

\author[BUW]{Matthias Ehrhardt\corref{cor}}
\cortext[cor]{Corresponding author}
\ead{ehrhardt@uni-wuppertal.de}

\author[UMI]{Moulay Rchid Sidi Ammi}
\ead{m.sidiammi@umi.ac.ma}

\address[BUW]{Applied and Computational Mathematics, University of Wuppertal,\\
  Gaußstraße 20, 42119 Wuppertal, Germany}
\address[UMI]{AMNEA Team, Department of Mathematics, FST Errachidia,\\
  Moulay Ismail University of Meknes, 52000 Errachidia, Morocco}

\begin{abstract}
We propose a conservative, semi-implicit nonstandard finite difference scheme for the reaction--diffusion equation $u_t=\lambda(|u_x|^{p-2}u_x)_x+f(u)$, $p>1$, where, as is customary, the $p$-Laplacian is regularized by a parameter $\eps>0$.
The scheme combines a denominator function, a lagged diffusion coefficient, and a nonlocal treatment of the loss part of the reaction.
For every time step and mesh size we prove that the scheme is well posed and positive, that it preserves an invariant region $[0,K]$ whenever the reaction admits one, and that it reproduces every constant equilibrium exactly.
Without reaction, it conserves mass, satisfies a discrete maximum principle, and, for $1<p\le 2$, decreases a discrete energy. 
An example shows that the energy may increase for $p>2$.
The truncation error is $O(k+h^2)$ at interior nodes.
Numerical tests confirm first order in time and second order in space, and show that the scheme remains positive and bounded for step sizes at which the explicit method breaks down.
\end{abstract}


\end{frontmatter}


\section{Introduction}\label{sec:intro}
Reaction--diffusion equations with $p$-Laplacian-type diffusion arise in models of non-Newtonian flow, filtration, and nonlinear heat conduction \cite{dibenedetto1993, Vazquez2006, Barrett1993, Zinihi2025FDE}.
Since the unknown is usually a density or a concentration, a numerical method should return nonnegative values and, if the reaction has a carrying capacity, respect it.
Explicit finite differences cannot guarantee this.
When the explicit scheme is written in conservative form, the coefficient of $u_i^m$ in the update is $1-\lambda k\,(c_{i+1/2}+c_{i-1/2})/h^2$, where the conductance $c\approx|u_x|^{p-2}$ is evaluated at the old time level.
Positivity can therefore be guaranteed only under a step size restriction that depends on the gradients, which are not known in advance, and this restriction becomes severe as $p$ grows (Section~\ref{sec:numerics}).

\textit{Nonstandard finite difference} (NSFD) methods address such failures with two devices:
a denominator function $\varphi(k)=k+O(k^2)$ that replaces $k$ in the difference quotient, and nonlocal (two-level) approximations of nonlinear terms 
\cite{Mickens1999, Mickens2006, Zinihi2026pLaplacianNSFD, Ehrhardt2013}.
Evaluating the diffusion coefficient at the old level (``lagged diffusivity'') is a classical device in image processing and in the iterative solution of $p$-Poisson problems \cite{Chan1999, Diening2020}.
Here we combine these ingredients:
lagging the conductance and treating the loss part of the reaction nonlocally yields a scheme whose system matrix is a tridiagonal $M$-matrix.
Unconditional positivity, an invariant region, exact equilibria, and an energy inequality then follow from short arguments, for every $p>1$ and for a whole class of reactions.
We also show where the argument stops: for $p>2$ the discrete energy can increase.

We do not claim that the individual devices are new.
In the previous work~\cite{Zinihi2026pLaplacianNSFD} we proposed an explicit NSFD scheme for the same class of problems, in which only the central node of each discrete flux is taken at the new time level.
That scheme is positive for every step size, but its spatial operator is not in conservative form, it requires $u^m>0$, the conductance is not regularized for $p<2$, and the $\ell^\infty$ bound proved there grows with the number of time steps.
The present scheme evaluates the conductance at cell midpoints in conservative form and lags it in time, so that the entire new time level enters through a tridiagonal $M$-matrix.
At the cost of one tridiagonal solve per step, this yields properties that the explicit scheme does not provide:
invariance of $[0,K]$ uniformly in $k$, $h$ and $m$, mass conservation and a discrete maximum principle, the exact logistic flow in the reaction part (Remark~\ref{rem:logistic}), and a discrete energy inequality for $1<p\le2$ together with an example showing that it fails for $p>2$.
We also give numerical evidence, including a comparison with iterating the lagged conductance to convergence (Table~\ref{tab:picard}).


\section{The Scheme}\label{sec:scheme}

For $x\in(0,1)$ and $t>0$ we consider a \textit{reaction--diffusion equation} of the form
\begin{equation}\label{eq:pde}
  u_t=\lambda\,\bigl(F_\eps(u_x)\bigr)_x+f(u),
  \qquad
  F_\eps(s)=(s^2+\eps^2)^{\frac{p-2}{2}}s,
\end{equation}
with zero-flux boundary conditions $F_\eps(u_x)=0$ at $x=0,1$, initial data $u(\cdot,0)=u_0\ge 0$, $\lambda>0$, $p>1$, and a fixed regularization parameter $\eps>0$, as in finite element treatments of the $p$-Laplacian \cite{Barrett1993}.
For $p=2$, \eqref{eq:pde} reduces to the semilinear heat equation for every $\eps$.
The \textit{reaction} is split as
\begin{equation}\label{eq:split}
  f(s)=a(s)-s\,b(s), \qquad a,b\in C[0,\infty),\quad a,b\ge 0.
\end{equation}
For the logistic term $f(s)=s(1-s)$ one takes $a(s)=b(s)=s$;
for $f\in C^1[0,\infty)$ with $f(0)\ge0$ one may take $a=f^+$ and $b(s)=f^-(s)/s$, extended continuously to $s=0$.

Let $h=1/N$, $x_i=ih$, $D_iv=(v_{i+1}-v_i)/h$, and let $w_0=w_N=h/2$ and $w_i=h$ for $0<i<N$.
For a grid function $u$ we define the 
\textit{conductances} and the \textit{discrete diffusion operator}
\begin{equation*}
  c_{i+1/2}[u] = \bigl( (D_iu)^2 +\eps^2 \bigr)^{\frac{p-2}{2}}>0,
  \qquad
  (L[u]v)_i = \frac{1}{w_i} \bigl( c_{i+1/2}D_i v-  c_{i-1/2} D_{i-1}v \bigr),
\end{equation*}
with the zero-flux convention $c_{-1/2}D_{-1}v=c_{N+1/2}D_Nv=0$.
With the denominator function $\varphi(k)=(\mathrm{e}^{\mathcal{L}k}-1)/\mathcal{L}$, $\mathcal{L}>0$, which satisfies $\varphi(k)=k+O(k^2)$ \cite{Mickens2006}, the \textit{NSFD scheme} reads
\begin{equation}\label{eq:scheme}
  \frac{u_i^{m+1}-u_i^m}{\varphi(k)}
  =\lambda\,\bigl( L[u^m] u^{m+1} \bigr)_i + a(u_i^m) - b(u_i^m)\,u_i^{m+1},
  \qquad i=0,\dots,N.
\end{equation}
The conductance is evaluated at the old level, the solution enters linearly at the new level, and the loss term $s\,b(s)$ is replaced by $b(u^m)u^{m+1}$.
Writing $\varphi=\varphi(k)$ and $\theta_i=\varphi\lambda/(w_ih)$, we can rewrite~\eqref{eq:scheme} as $B^mu^{m+1}=u^m+\varphi\,a(u^m)$ with the tridiagonal matrix
\begin{equation*}
  B^m_{ii} = 1+\varphi\, b(u_i^m)+\theta_i\,(c_{i+1/2}+c_{i-1/2}), \qquad
  B^m_{i,i+1} = -\theta_i\,c_{i+1/2}, \qquad
  B^m_{i,i-1} = -\theta_i\,c_{i-1/2},
\end{equation*}
where $c_{-1/2}=c_{N+1/2}=0$.
Hence, one step costs $O(N)$ operations with the Thomas algorithm 
\cite{thomas1949elliptic} 
the same order as an explicit step.

\section{Properties of the Scheme}\label{sec:properties}
We now collect the main properties of scheme~\eqref{eq:scheme} that hold for every time step and mesh size.
\begin{theorem}\label{thm:main}
Let $p>1$ and $\eps,\lambda,k,h>0$ be arbitrary, let the splitting of the reaction~\eqref{eq:split} hold, and let $u^m\ge0$ be a grid function.
\begin{enumerate}[(a)]
\item $B^m$ is a strictly diagonally dominant $M$-matrix with $\|(B^m)^{-1}\|_\infty \le 1$. Hence~\eqref{eq:scheme} has a unique solution $u^{m+1}$.
\item $u^{m+1}\ge0$, and $u^{m+1}>0$ if $u^m>0$.
\item If $a(s)\le K\,b(s)$ on $[0,K]$ for some $K>0$ and $0\le u^0\le K$, then $0\le u^m\le K$ for all $m$.
\item If $f(u^*)=0$ for some $u^*\ge0$ and $u^m\equiv u^*$, then $u^{m+1}\equiv u^*$.
\item If $f\equiv0$, then $\sum_i w_i u_i^m$ is independent of $m$, and $\min_iu_i^0\le u_i^m \le \max_i u_i^0$.
\end{enumerate}
\end{theorem}
\begin{proof}
(a) $B^m$ has a positive diagonal, non-positive off-diagonal entries, and $B^m_{ii}-\sum_{j\ne i}|B^m_{ij}|=1+\varphi\, b(u_i^m)\ge1$.
For any vector $v$, choose $i$ with $|v_i|=\|v\|_\infty$.
Then
$|(B^mv)_i| \ge B^m_{ii}|v_i| - \sum_{j\ne i}|B^m_{ij}|\,|v_j| \ge \bigl( B^m_{ii} - \sum_{j\ne i}|B^m_{ij}| \bigr)\|v\|_\infty \ge \|v\|_\infty$, so $\|B^mv\|_\infty \ge \|v\|_\infty$.
Hence $B^m$ is nonsingular and $\|(B^m)^{-1}\|_\infty \le 1$.
Writing $B^m=D-R$ with $D$ diagonal and $R\ge0$, diagonal dominance gives $\|D^{-1}R\|_\infty<1$, 
so $(B^m)^{-1} = \sum_{j\ge0} (D^{-1}R)^j D^{-1} \ge 0$ entry-wise, with diagonal entries at least $1/D_{ii}$.

(b) Since $a\ge0$ and $\varphi>0$, the right-hand side satisfies $u^m+\varphi\,a(u^m)\ge u^m\ge0$.
As $(B^m)^{-1}\ge0$, we get $u^{m+1}\ge0$.
If $u^m>0$, then $u_i^{m+1}\ge\bigl(u_i^m+\varphi\, a(u_i^m)\bigr)/D_{ii}>0$.

(c) We argue by induction. To do so, assume $0\le u^m\le K$.
By (b), $u^{m+1}\ge0$.
Let $M\coloneqq\max_iu_i^{m+1}=u^{m+1}_{i^*}$ and $\beta\coloneqq b(u_{i^*}^m)$.
In row $i^*$ of~\eqref{eq:scheme} the diffusion part $\theta_{i^*}\sum_{j=i^*\pm1} c_{(i^*+j)/2}\,(u_{i^*}^{m+1}-u_j^{m+1})$ is nonnegative, so
$(1+\varphi\beta)M\le u_{i^*}^m+\varphi\, a(u_{i^*}^m)\le K+\varphi K\beta$,
where we used $u_{i^*}^m\le K$ and $a\le Kb$ on $[0,K]$. Hence $M\le K$.

(d) Constants lie in the kernel of the diffusion part, and $a(u^*)=u^*b(u^*)$.
Hence $B^mu^*=(1+\varphi\, b(u^*))u^*=u^*+\varphi\, a(u^*)$, i.e., $u^*$ solves~\eqref{eq:scheme}.

(e) With $a=b=0$, multiplying~\eqref{eq:scheme} by $w_i$ and summing over $i$, the diffusion terms telescope to zero by the zero-flux convention.
At a maximizer $i^*$ of $u^{m+1}$, row $i^*$ gives $u^{m+1}_{i^*}\le u^m_{i^*}\le\max_i u_i^m$. 
The minimum is treated in the same way.
\end{proof}
\begin{remark}\label{rem:conductance}
The proof of Theorem~\ref{thm:main} uses only the positivity of the conductances $c_{i+1/2}$, not the way they are computed.
Parts (a)--(e) therefore hold for any positive conductance in place of $c_{i+1/2}[u^m]$, in particular for every iterate of a Picard iteration that updates the conductance within a time step.
Theorem~\ref{thm:energy} below, in contrast, does use $c=c[u^m]$.
\end{remark}

\begin{remark}\label{rem:logistic}
For the logistic reaction with $a(s)=b(s)=s$ and $\mathcal{L}=1$, the reaction part of~\eqref{eq:scheme} is the map $u\mapsto u(1+\varphi)/(1+\varphi u)=u\,\mathrm{e}^k/\bigl(1+u(\mathrm{e}^k-1)\bigr)$, which is the exact flow of $u'=u(1-u)$.
Part (c) holds with $K=1$.
\end{remark}

For the energy estimate, let $G(s)=(s^2+\eps^2)^{p/2}/p$, so that $G'=F_\eps$, and define
\begin{equation*}
  E_h(v) = h \sum_{i=0}^{N-1} G(D_iv), \qquad
  \|v\|_W^2 = \sum_{i=0}^N w_i v_i^2 .
\end{equation*}

\begin{theorem}\label{thm:energy}
Let $f\equiv0$ and $1<p\le2$. Then, for every $k,h,\eps,\lambda>0$,
\begin{equation*}
  E_h(u^{m+1}) + \frac{1}{\lambda\varphi(k)}\|u^{m+1}-u^m\|_W^2\le E_h(u^m).
\end{equation*}
\end{theorem}

\begin{proof}
Write $u^+=u^{m+1}$, $u^-=u^m$, $z=u^+-u^-$, and $c_i=c_{i+1/2}[u^-]$.
Summation by parts with the zero-flux convention gives $\sum_i w_i z_i (L[u^-]u^+)_i = -h \sum_ic_i\,D_iz\,D_iu^+$.
Multiplying~\eqref{eq:scheme} by $w_iz_i$ and summing yields
\begin{equation}\label{eq:test}
  h \sum_i c_i\,D_iz\,D_iu^+ = -\frac{1}{\lambda\varphi}\|z\|_W^2 .
\end{equation}
Let $g(t)=(t+\eps^2)^{p/2}/p$, so that $G(s)=g(s^2)$.
For $p\le2$ the function $g$ is concave, hence $G(s')\le G(s)+g'(s^2)(s'^2-s^2)$ with $g'(s^2) = \tfrac{1}{2}(s^2+\eps^2)^{(p-2)/2}$.
Applying this with $s=D_iu^-$ and $s'=D_iu^+$ gives $E_h(u^+) \le E_h(u^-) + \tfrac{h}{2} \sum_i c_i \bigl( (D_iu^+)^2-(D_iu^-)^2 \bigr)$.
Since $(D_iu^+)^2-(D_iu^-)^2 = 2D_iz\,D_iu^+-(D_iz)^2$, we obtain $E_h(u^+) \le E_h(u^-) + h \sum_ic_i\,D_iz\,D_iu^+-\tfrac{h}{2} \sum_i c_i (D_iz)^2$, and the claim follows from~\eqref{eq:test} after dropping the last (non-positive) term.
\end{proof}

\begin{remark}\label{rem:energy-fails}
For $p>2$ the function $g$ is convex, and not only the argument but also the conclusion fails.
Take $p=3$, $N=40$, $\eps=0.01$, $\lambda=1$, $k=0.1$, and $u_0(x) = 1 + \tfrac{1}{2} \cos 7\pi x + \tfrac{1}{2} \mathrm{e}^{-60(x-1/2)^2}$.
In the third step, $E_h$ rises by about $22\,\%$, from about $1.7\cdot10^{-4}$ to $2.1\cdot10^{-4}$ (Figure~\ref{fig:energy}, right).
\end{remark}

\begin{proposition}\label{prop:consistency}
Let $a$ and $b$ be Lipschitz continuous, and let $u$ be a solution of~\eqref{eq:pde} that is $C^4$ in $x$ and $C^2$ in $t$.
Inserting $u$ into~\eqref{eq:scheme} leaves a residual $\tau_i^m=O(k+h^2)$ at interior nodes and $\tau_i^m=O(k+h)$ for $i\in\{0,N\}$, with constants depending on $p$, $\eps$, $\lambda$, and norms of $u$.
\end{proposition}
\begin{proof}
Since $\varphi(k)=k+O(k^2)$, the time quotient equals $u_t+O(k)$.
Because $c_{i+1/2}\le\eps^{p-2}$ for $p<2$ and $c_{i+1/2}$ is bounded on bounded gradients for $p\ge2$, replacing $u(t_{m+1})$ by $u(t_m)$ in the diffusion and reaction terms costs $O(k)$.
For the spatial part, $D_iu=u_x(x_{i+1/2})+\tfrac{h^2}{24}u_{xxx}(x_{i+1/2})+O(h^3)$ and $F_\eps$ is smooth, so $c_{i+1/2}D_iu=F_\eps(u_x)(x_{i+1/2})+h^2\psi(x_{i+1/2})+O(h^3)$ with $\psi\in C^1$.
The difference quotient of this flux is $(F_\eps(u_x))_x+O(h^2)$ at interior nodes.
At $i=0$ the half cell gives $(2/h)F_\eps(u_x)(h/2)+O(h)=(F_\eps(u_x))_x(0)+O(h)$, using $F_\eps(u_x(0))=0$.
The case $i=N$ is analogous.
\end{proof}

\section{Numerical Experiments}\label{sec:numerics}
Unless stated otherwise, all runs use $\mathcal{L}=1$, $\lambda=1$, and the logistic reaction with $a(s)=b(s)=s$.
%
We use the manufactured solution $u(x,t) = 1 + \tfrac{1}{2} \mathrm{e}^{-t}\cos\pi x$ and add to~\eqref{eq:pde} the source term $g$ that makes it exact for the regularized operator, $g^+$ is included in $a$ and $g^-/u$ in $b$.
Tables~\ref{tab:time} and~\ref{tab:space} show first order in time and second order in space for $p=1.5$ and $p=3$.
The last spatial rates fall below two because the first-order time error is still visible: for $p=3$ and $N=320$, halving $k$ changes the error from $2.85\cdot10^{-6}$ to $2.55\cdot10^{-6}$.

\begin{table}[!htb]
\centering\footnotesize
\setlength{\tabcolsep}{0.25cm}
\begin{minipage}[t]{0.48\textwidth}\vspace{0pt}\centering
\captionof{table}{Time refinement with $N=320$: maximum-norm errors and observed orders ($\eps=0.2$, $T=0.2$).}\label{tab:time}
\begin{tabular}{lcccc}
\toprule
 & \multicolumn{2}{c}{$p=1.5$} & \multicolumn{2}{c}{$p=3$}\\
\cmidrule(lr){2-3}\cmidrule(lr){4-5}
$k$ & Error & Order & Error & Order\\
\midrule
$T/2^{1}$ & $5.10\times10^{-2}$ & --   & $3.88\times10^{-2}$ & --\\
$T/2^{2}$ & $3.25\times10^{-2}$ & 0.65 & $2.02\times10^{-2}$ & 0.94\\
$T/2^{3}$ & $1.86\times10^{-2}$ & 0.80 & $1.05\times10^{-2}$ & 0.94\\
$T/2^{4}$ & $1.00\times10^{-2}$ & 0.90 & $5.37\times10^{-3}$ & 0.97\\
$T/2^{5}$ & $5.19\times10^{-3}$ & 0.95 & $2.72\times10^{-3}$ & 0.98\\
$T/2^{6}$ & $2.64\times10^{-3}$ & 0.97 & $1.36\times10^{-3}$ & 0.99\\
$T/2^{7}$ & $1.33\times10^{-3}$ & 0.98 & $6.84\times10^{-4}$ & 1.00\\
\bottomrule
\end{tabular}
\end{minipage}
\hfill
\begin{minipage}[t]{0.48\textwidth}\vspace{0pt}\centering
\captionof{table}{Space refinement with $k=T/2^{17}$: maximum-norm errors and observed orders ($\eps=0.2$, $T=0.2$).}\label{tab:space}
\begin{tabular}{lcccc}
\toprule
 & \multicolumn{2}{c}{$p=1.5$} & \multicolumn{2}{c}{$p=3$}\\
\cmidrule(lr){2-3}\cmidrule(lr){4-5}
$N$ & Error & Order & Error & Order\\
\midrule
10  & $3.95\times10^{-2}$ & --   & $3.69\times10^{-3}$ & --\\
20  & $3.18\times10^{-3}$ & 3.63 & $5.57\times10^{-4}$ & 2.73\\
40  & $2.87\times10^{-4}$ & 3.47 & $1.45\times10^{-4}$ & 1.94\\
80  & $6.99\times10^{-5}$ & 2.04 & $3.67\times10^{-5}$ & 1.98\\
160 & $1.84\times10^{-5}$ & 1.92 & $9.63\times10^{-6}$ & 1.93\\
320 & $5.59\times10^{-6}$ & 1.72 & $2.85\times10^{-6}$ & 1.76\\
\bottomrule
\end{tabular}
\end{minipage}
\end{table}

Table~\ref{tab:explicit} compares scheme~\eqref{eq:scheme} with the explicit method based on the same conservative spatial discretization.
The initial datum is a narrow bump on a nearly vanishing background, and $h=0.01$.
For $p=3$ the explicit method blows up even for $k=10^{-4}$, since the conductance grows with the gradient.
For $p=1.5$ it produces negative values for every $k$ we tried.
For $k=10^{-4}$ the violation is small ($\min u=-3.2\cdot10^{-5}$ at intermediate times, and the solution is positive again at $T$), but it becomes large for larger $k$.
Scheme~\eqref{eq:scheme} stays in $[0,1]$ in all cases, as predicted by Theorem~\ref{thm:main}.
Note that for large $k$ we claim only positivity and boundedness, the accuracy is still $O(k)$.

\begin{table}[!htb]
\centering\footnotesize
\setlength{\tabcolsep}{0.35cm}
\caption{Explicit FDM versus scheme~\eqref{eq:scheme} for $u_0(x)=0.02+0.98\,\mathrm{e}^{-200(x-1/2)^2}$, $\eps=0.05$, $h=0.01$, $T=0.05$; ``blow-up'' means $\|u\|_\infty>10^4$.}\label{tab:explicit}
\begin{tabular}{llll}
\toprule
$p$ & $k$ & Explicit FDM & Scheme~\eqref{eq:scheme}: $[\min u,\max u]$ at $T$\\
\midrule
1.5 & $10^{-4}$ & negative values ($\min u=-3.2\times10^{-5}$) & $[0.1388, 0.1582]$\\
1.5 & $10^{-3}$ & negative values ($\min u=-2.1$)             & $[0.1369, 0.1603]$\\
1.5 & $10^{-2}$ & negative values ($\min u=-311.9$)           & $[0.1261, 0.1832]$\\
3   & $10^{-4}$ & blow-up at $t=0.0003$                       & $[0.0764, 0.2211]$\\
3   & $10^{-3}$ & blow-up at $t=0.0030$                       & $[0.0753, 0.2222]$\\
3   & $10^{-2}$ & blow-up at $t=0.0200$                       & $[0.0544, 0.2432]$\\
\bottomrule
\end{tabular}
\end{table}

The lag can be removed by a Picard iteration within each step: starting from $v^0=u^m$, the iterate $v^{j+1}$ solves~\eqref{eq:scheme} with $c[u^m]$ replaced by $c[v^j]$, until convergence.
The reaction is treated as in~\eqref{eq:scheme}.
By Remark~\ref{rem:conductance} every iterate is positive, so positivity does not distinguish the two approaches --  the difference lies in the cost.
Table~\ref{tab:picard} reports the mean number of tridiagonal solves per step needed to reach $\|v^{j+1}-v^j\|_\infty<10^{-8}$ for the data of Table~\ref{tab:explicit}.
For $p=1.5$ it ranges from about 5 to 17: for $p=3$ it grows quickly with $k$, and the cap of 200 iterations was reached in some steps.
The solution of~\eqref{eq:scheme} differs from the converged one by an amount of order $k$, as expected from lagging the conductance.
In this sense, \eqref{eq:scheme} is the one-iteration member of this family.

\begin{table}[!htb]
\centering\footnotesize
\setlength{\tabcolsep}{0.6cm}
\caption{Scheme~\eqref{eq:scheme} (solution $u$) versus the converged Picard iteration (solution $u_{\mathrm{Pic}}$) within each step, stopped when $\|v^{j+1}-v^{j}\|_\infty<10^{-8}$ or after 200 iterations; data as in Table~\ref{tab:explicit}. 
Entries marked $\ge$ are lower bounds, since the cap was reached in some steps.}\label{tab:picard}
\begin{tabular}{llcc}
\toprule
$p$ & $k$ & Picard solves per step & $\|u-u_{\mathrm{Pic}}\|_\infty$ at $T$\\
\midrule
1.5 & $10^{-4}$ & 5.4         & $1.5\times10^{-4}$\\
1.5 & $10^{-3}$ & 10.1        & $1.4\times10^{-3}$\\
1.5 & $10^{-2}$ & 17.2        & $1.5\times10^{-2}$\\
3   & $10^{-4}$ & 4.6         & $1.9\times10^{-4}$\\
3   & $10^{-3}$ & $\ge 28.9$  & $6.8\times10^{-4}$\\
3   & $10^{-2}$ & $\ge 168$   & $1.2\times10^{-2}$\\
\bottomrule
\end{tabular}
\end{table}

Table~\ref{tab:reaction} isolates the effect of the reaction term.
We keep the denominator function and the lagged implicit diffusion of~\eqref{eq:scheme} but treat the reaction explicitly.
For data above $K=1$ and a small diffusion coefficient, this variant produces negative values for $k=0.5$ ($\min u=-0.20$) and reaches $-181$ for $k=1$.
The nonlocal treatment keeps $u\ge1$ throughout, with final state in $[1,1.0313]$ for all $k$ tested.
This lower bound follows from the minimum version of the argument in Theorem~\ref{thm:main}(c), since $a=b$.

\begin{table}[!htb]
\centering\footnotesize
\setlength{\tabcolsep}{0.6cm}
\caption{Smallest value of $u$ over the whole run for data above $K=1$: $u_0(x)=1+2\,\mathrm{e}^{-50(x-1/2)^2}$, $p=3$, $\lambda=10^{-3}$, $\eps=0.05$, $N=100$, $T=3$. Both columns use the same $\varphi(k)$ and the same lagged diffusion.}\label{tab:reaction}
\begin{tabular}{lcc}
\toprule
$k$ & Explicit reaction & Nonlocal reaction~\eqref{eq:scheme}\\
\midrule
0.1  & $1.0000$    & $1.0000$\\
0.25 & $1.0000$    & $1.0000$\\
0.5  & $-0.2017$   & $1.0000$\\
1    & $-181.1501$ & $1.0000$\\
\bottomrule
\end{tabular}
\end{table}

\begin{figure}[!htb]
\centering
\includegraphics[width=0.9\textwidth]{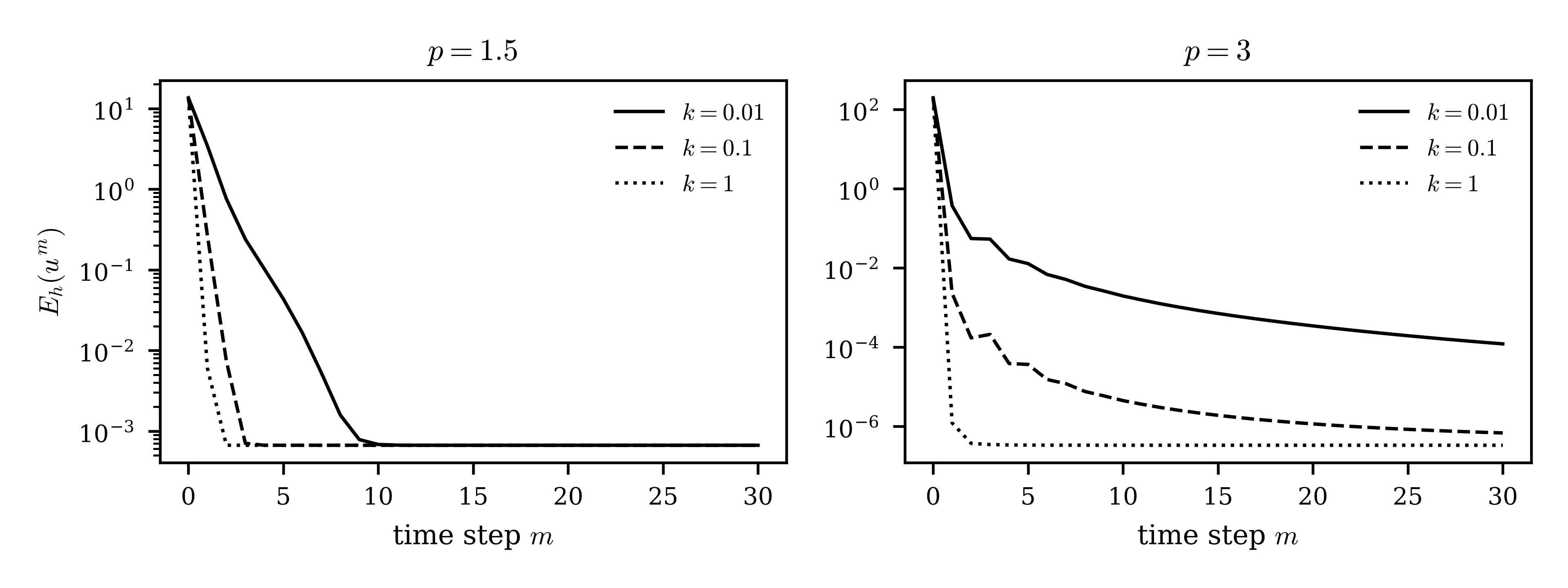}
\caption{Discrete energy $E_h(u^m)$ for $f=0$, $N=40$, $\eps=0.01$, $\lambda=1$, and $u_0(x)=1+\tfrac{1}{2} \cos 7\pi x+\tfrac{1}{2}\mathrm{e}^{-60(x-1/2)^2}$.
Left: $p=1.5$ (monotone for all $k$). Right: $p=3$ (an increase occurs for $k=0.1$).}\label{fig:energy}
\end{figure}

Figure~\ref{fig:energy} shows the discrete energy $E_h(u^m)$ for $f=0$ and the initial datum of Remark~\ref{rem:energy-fails}.
For $p=1.5$ the energy decreases monotonically for $k=0.01,0.1,1$, as guaranteed by Theorem~\ref{thm:energy}.
For $p=3$ the same data give a nonmonotone history for $k=0.1$ (Remark~\ref{rem:energy-fails}).
In randomized tests with $3000$ steps for each $p$ (random data, $k\in[10^{-4},10]$, $\eps\in[10^{-3},10^{-0.5}]$), the energy never increased for $p\in\{1.5,2,2.5\}$, but it increased in $95$ steps for $p=3$ and in $250$ steps for $p=5$.

\section{Conclusion}\label{sec:conclusion}
Scheme~\eqref{eq:scheme} guarantees positivity, an invariant region, exact constant equilibria, mass conservation, and a maximum principle for every step size, as well as an energy inequality for $1<p\le2$, at the cost of one tridiagonal solve per step.
Several questions remain open.
Convergence of~\eqref{eq:scheme} to the solution of the regularized problem is supported by Proposition~\ref{prop:consistency} and Tables~\ref{tab:time}--\ref{tab:space}, but is not proved here.
The standard error analysis is hampered by the dependence of the conductance on the discrete gradient, which introduces a factor $h^{-1}$.
The dependence of the error constants on $\eps$ and the limit $\eps\to0$ are not analyzed.
The scheme is only first-order accurate in time, and since $\varphi(k)$ grows exponentially in $\mathcal{L}k$, $\mathcal{L}$ should be kept moderate.
In two space dimensions the $M$-matrix structure persists, but a sparse linear system has to be solved in each step.
Finally, an energy-stable variant for $p>2$ remains open.
Remark~\ref{rem:energy-fails} shows that it requires a modification of~\eqref{eq:scheme}.

\section*{Declarations}

\subsection*{Data availability}
No data was used for the research described in this letter.

\subsection*{Declaration of competing interest}
The authors declare that they have no known competing financial interests or personal relationships that could have appeared to influence the work reported in this letter.

\bibliographystyle{elsarticle-num}
\bibliography{References}

\end{document}